\documentclass[a4paper, 11pt]{amsart}
\usepackage{amssymb}
\usepackage{a4wide}
\usepackage{verbatim}
\usepackage{graphics}

\usepackage{epsfig}
\usepackage{MnSymbol}

\usepackage{enumerate}

\usepackage{hyperref}
\newtheorem{theorem}{Theorem}

\theoremstyle{definition}

\newtheorem{example}[theorem]{Example}

\newtheorem{remark}[theorem]{Remark}
\newtheorem{proposition}[theorem]{Proposition}
\newtheorem{definition}[theorem]{Definition}

\newcommand{\pd}[2]{\frac{\partial#1}{\partial#2}}
\newcommand{\vp}[1]{\frac{\partial}{\partial#1}}

\newcommand{\X}{\mathfrak{X}}

\newcommand{\Sp}{\text{Sp}}
\newcommand{\R}{\mathbb{R}}

\begin{document}
\title{Explicit solutions of PDE via Vessiot theory and solvable structures}
\author{Naghmana Tehseen}
\address{ Department of Mathematics and Statistics, La Trobe University, Victoria, 3086, Australia}
\email{n.tehseen@latrobe.edu.au,naghmanatehseen@gmail.com}
\date{\today}
\noindent \begin{abstract}
We consider the problem of computing the integrable sub-distributions of the non-integrable Vessiot distribution of multi-dimensional second order partial differential equations (PDEs). We use Vessiot theory and solvable structures to find the largest integrable distributions contained in the Vessiot distribution associated to second order PDEs. In particular, we show how the solvable symmetry structure of the original PDE can be used to construct integrable sub-distributions leading to group invariant solutions of the PDE in two and more than two independent variables.
\end{abstract}
\maketitle
\section{Introduction}
\noindent In a recent contribution \cite{NTGP13}, the authors extended the results of \cite{sherr} on the use of solvable structures for integrating ODEs (that is, one dimensional distributions), and the work of \cite{MbarcoGP01,MBGP01} on first order PDEs, to the integration of second order PDEs of one dependent variable and two independent variables. The equations are of the form
\begin{equation*}
  u_{yy}=F(x,y,u,u_x,u_y,u_{xx},u_{xy})
\end{equation*}
for some smooth function $F.$ Such equations may be formulated \cite{Vessiot1924} as a distribution of vector fields on a jet space $J^2(\R^2,\R)$ and this distribution is called a Vessiot distribution \cite{Fackerell85}. The maximal integrable sub-distributions of the  Vessiot distribution which satisfy an appropriate independence condition represent local solutions of the PDE system. The authors provided a systematic approach for solving second order linear and non-linear PDEs in the presence of symmetries, a general class of symmetries known as \emph{solvable structure}, which are not necessarily of point type. They showed  how to use a solvable  structure to integrate a PDE in the original coordinates, but most importantly, they showed how to impose a solvable structure (referred as ``group invariance method") on a PDE so as to determine particular largest solvable sub-distributions of the Vessiot distribution.

In this paper, we enlarge the scope of applicability of the preceding study by allowing more independent and dependent variables into the picture. We apply the same geometric technique to find the largest integrable sub-distributions of the Vessiot distribution of the multi-dimensional second order PDEs. We demonstrate the higher dimensional version of the group invariance method. The basic tools that we use are Vessiot's geometric formulation of differential equations \cite{Vessiot1924} and the integrating factor technique developed in \cite{sherr}.  We will indeed see that this technique can be applied to second order PDEs in any dimensions. The main difficulty arising here is the complexity of differential conditions, and the volume of symbolic calculations required. To overcome this problem we apply the group invariance method and show how one can construct the ``group invariant" solutions for PDEs in two and more than two independent variables. To generate such solvable structures for each integrable sub-distribution we use the symmetry determination software package DIMSYM \cite{sherrdimsym} operating as a REDUCE \cite{Reduce}  overlay. We also use the exterior calculus package EXCALC \cite{Sexcalc}.

The plan of the paper is to first provide some key concepts of Vessiot theory. We then describe a method to find the largest integrable sub-distributions which satisfy an appropriate independence condition. We also demonstrate the construction of group invariant solutions of multi-dimensional second order PDEs. We then apply the method to second order PDEs of one dependent variable and three independent variables. In the next part of the paper, we apply the method to second order PDEs of one dependent variable and four independent variables and to systems of PDEs. To illustrate the method some concrete examples are given.

\section{The PDE problem}
\noindent Consider a system of partial differential equations  of $m$ independent variables $x^i$ and $n$ dependent variables  $u^j$,
\begin{align}F^a(x^i,u^j,u^j_{i_1},u^j_{i{_1}i{_2}}\ldots u^j_{i_1\ldots i_k})=0. \label{eqn:g.pde}\end{align} Denote by $X$ the space of independent variables and by $U$ the space of dependent variables. The subscripts $1\leq i_1\leq\ldots \leq i_k\leq m$ are used to specify the partial derivatives of $u^j$, where $k$ is the maximum order of the system. For definitions and notations we refer to the monograph \cite{Stormark}.

\noindent Consider the trivial bundle $\pi:X\times U\rightarrow X$ and let $f$ and $g$ be two smooth sections of $\pi$. We say that $f$ and $g$  are equivalent to order $k$ at $x$ if and only if  \[\frac{\partial^{q_1+ \ldots +q_m} f}{(\partial x^1)^{q_1}\ldots(\partial x^m)^{q_m}}(x)=\frac{\partial^{q_1+\ldots+q_m} g}{(\partial x^1)^{q_1}\ldots(\partial x^m)^{q_m}}(x)\] for all $m-$tuples $(q_1\ldots q_m)$ with $q_1+\ldots+q_m\leq k.$
The equivalence class of a smooth section $f$ at a point $x$ is the $k-$jet of $f$ at $x$  and is denoted by $j^k_x f.$ The set of all $k-$jets of smooth sections of $\pi$ constitutes the bundle of $k-$jets of maps $X\rightarrow U$ and denoted by $J^k(X,U).$ The zeroth order jet bundle, $J^0(X,U),$ is identified with $X\times U.$ The $k-$graph of smooth section $f$ of $\pi$ is the map $j^k f: X\rightarrow J^k(X,U)$ defined by $x\rightarrow j^k_x f.$ Thus, the image of $k-$graph of a section is a $m-$dimensional immersed submanifold of $J^k(X,U).$

The geometry of jet bundles is to a large extent determined by their \emph{contact structure.} The contact co-distribution consists of all one$-$forms such that their pull-back by a prolonged section vanishes. Locally, it is spanned by the contact forms
\begin{align}\label{contactsystem}
\theta^j_I:=d u^j_I-\sum_{i=1}^m u_{I, i}^j d x^i,
\end{align} where $I$ is a multi-index of order less than or equal to $k-1$ and is denoted by $\Omega^k(X,U).$

The dual distribution $\Omega^k(X,U)^\perp$ consisting of all vector fields annihilated by contact co-distribution $\Omega^k(X,U)$. A straightforward calculation \cite{Fesser09} shows that it is generated by
\begin{align}\label{Vessiot dist}
V_{i}^{(k)}&:=\partial_i+\sum_{j=1}^n \sum_{0\leq \left|I\right|<k}u_{I,i}^j\partial_{u^j_I},\quad 1\leq i\leq m,\nonumber \\
V_{j}^{I}&:=\partial_{u^j_I},\quad \quad \left|I\right|=k,\quad \quad 1\leq j \leq n.
\end{align}
The Vessiot distribution is the restriction of the above contact distribution to the submanifold defined by our differential condition $F = 0$.  If a distribution is Frobenius integrable then the integral manifold of the distribution can be found by various constructive methods. In general, the Vessiot distribution is not Frobenius integrable.

Vessiot \cite{Vessiot1924} has given an algorithm for constructing all the Frobenius integrable sub-distributions of any given distribution. If the distribution is not Frobenius integrable, then it is interesting to ask for the sub-distributions which are Frobenius integrable. He looks for generic sub-distributions that satisfy the algebraic constraint called by him \emph{involutions of degree r} and then shows that maximal such involutions can be deformed so as to be Frobenius integrable. He does this via the Cauchy-Kowalevski theorem. We will give an equivalent method to find the largest integrable Vessiot sub-distributions. We choose to work with differential forms rather than vector fields as the integrating factor technique \cite{sherr} to integrate the integrable Vessiot sub-distributions is phrased this way.

 So the first problem is to locate the largest integrable sub-distributions which satisfy the transverse condition.  Then we apply the integrating factor technique \cite{sherr} to integrate these integrable Vessiot sub-distributions. We also show how solvable symmetry structure can be used to find the integrable sub-distributions.

 The main result of \cite{sherr} is reproduced below:
\begin{theorem} \label{sherr result}
 Let $\Omega$ be a decomposable $k-$form on a manifold $M^n,$ and let $\Sp(\{X_1,\ldots, X_k\})$ be a $k-$dimensional distribution on an open $U\subseteq M^n$ satisfying $X_i\righthalfcup \Omega\neq 0$ everywhere on $U.$ Further suppose that $\Sp(\{X_{j+1},\ldots, X_k\}\bigcup\text{ker}~\Omega)$ is integrable for some $j<k$ and that $X_i$ is a symmetry of $\Sp(\{X_{i+1},\ldots, X_{k}\}\bigcup \text{ker}~\Omega)~\text{for}~ i=1,\ldots, j.$

 Put $\sigma^i:=X_1\righthalfcup \ldots \righthalfcup \bar{X}_i\righthalfcup\ldots\righthalfcup X_k\righthalfcup\Omega,$
where $\bar{X}_i$ indicates that this argument is missing and $\omega^i:=\frac{\sigma^i}{X_i\righthalfcup \sigma^i}~ \text{for}~i=1,\ldots,k$
so that $\{\omega^1,\ldots, \omega^k\}$ is dual to $\{X_1,\ldots, X_k\}.$ Then $d \omega^1=0;~d \omega^2=0~\mod~\omega^1;~d \omega^3=0~\mod~\omega^1,\omega^2;~~~\ldots~;~ d \omega^j=0~\mod~\omega^1,\ldots,\omega^{j-1},$ so that locally
 \begin{align*}
   \omega^1&=d \gamma^1,\\
   \omega^2&=d \gamma^2-X_1(\gamma^2)d\gamma^1,\\
   \omega^3&=d \gamma^3-X_2(\gamma^3)d\gamma^2-(X_1(\gamma^3)
   -X_2(\gamma^3)X_1(\gamma^2))d\gamma^1,\\ \vdots\\
  \omega^{j}&\equiv d \gamma^{j}~\mod~d\gamma^1,\ldots, d\gamma^{j-1},
 \end{align*}
 for some $\gamma^1,\ldots,\gamma^j\in \bigwedge^0 T^* U.$ The system $\{\omega^{j+1},\ldots,\omega^k\}$ is integrable modulo $d\gamma^1,\ldots, d\gamma^j$ and locally $\Omega=\gamma^0 d\gamma^1\wedge d\gamma^2 \wedge \ldots \wedge d\gamma^j\wedge \omega^{j+1}\wedge\ldots \wedge \omega^k ~\text{for some}~\gamma^0\in \bigwedge^0(T^* U).$ Each $\gamma^i$ is uniquely defined up to the addition of arbitrary function of $\gamma^1,\ldots,\gamma^{i-1}.$
 \end{theorem}
\noindent We will now explain the Frobenius integrability of $p-$forms.
\begin{definition}
A differential $p-$form $\bar\Omega\in\bigwedge^p(M^n)$ is simple or decomposable if it is the wedge product of $p$ $1-$forms.
\end{definition}
\begin{definition}
A constraint $1-$form $\theta$ for differential form $\bar\Omega$ is a $1-$form satisfying $\theta\wedge{\bar\Omega}=0,$ which implies $Y\righthalfcup\theta=0,~ \forall Y\in \text{ker}~{\bar\Omega}.$
\end{definition}
\begin{definition}
A characterizing form for $m-$dimensional distribution $D$ is a form on $M^n$ of degree $(n-m)$ which is the exterior product of $(n-m)$ constraint forms.
\end{definition}
\begin{definition}
Let $\bar\Omega\in\bigwedge^p(M^n)$ for some $p>1$ be decomposable, $\bar\Omega$ is Frobenius integrable if $d\bar\Omega=\lambda\wedge\bar\Omega.$ Equivalently, $D:=\text{ker}\ \bar\Omega$ is Frobenius integrable. Note that $\text{dim(ker}\ \bar\Omega)= n-p,$ since $\bar\Omega$ is simple.
\end{definition}

\begin{definition} \label{def:sym cond}
   Let $\bar{\Omega}$ be a simple $p-$form on $M^n.$ Then a set of $p$ linearly independent vector fields $X_1,\dots,X_p\in \mathfrak{X}(M^n)$ forms a solvable symmetry structure for $\bar{\Omega}$ if the sequence of simple forms  $$\bar\Omega,\ X_1 \righthalfcup \bar\Omega,\ \dots\ ,\ X_{p-1}\righthalfcup \dots\righthalfcup  X_1\righthalfcup  \bar\Omega$$ satisfies
 \begin{align}\label{sym cond}
 \mathcal{L}_{X_1}\bar\Omega &=\ell_1\bar\Omega,\nonumber\\
 \mathcal{L}_{X_2}(X_1\righthalfcup\bar\Omega) &=\ell_2(X_1\righthalfcup\bar\Omega),\\
 \vdots\nonumber\\
 \mathcal{L}_{X_{p}}(X_{p-1}\righthalfcup\dots\righthalfcup X_1\righthalfcup\bar\Omega) &=\ell_{p}(X_{p-1}\righthalfcup\dots\righthalfcup X_1\righthalfcup\bar\Omega),\nonumber
 \end{align}
  for some smooth functions $\ell_1,\ldots,\ell_p.$
\end{definition}

  In \cite{sherr}, authors have extended the Lie's approach to integrate a Frobenius integrable distribution via a solvable structure of symmetries. In that paper, a Frobenius integrable distribution is given first. In our work, we start with non-integrable distribution and then we impose the symmetry conditions \eqref{sym cond} to construct a Frobenius integrable distribution.

Note that in definition \ref{def:sym cond}, we have not assumed the Frobenius integrability of $\bar{\Omega}.$ The way in which the symmetry conditions \eqref{sym cond} act as a catalyst  on $\bar{\Omega}$ can be seen more comprehensively by recalling the sequence \[\bar{\Omega},\ X_1 \righthalfcup \bar{\Omega},\ \dots\ ,\ X_{p-1}\righthalfcup\dots\righthalfcup X_1\righthalfcup \bar{\Omega}\]
are all simple forms and that they are all Frobenius integrable if $\bar\Omega$ is Frobenius integrable. Hence, the Frobenius integrability of each form in the sequence is a necessary condition for the Frobenius integrability of $\bar\Omega$. We can summarize this result as:
\begin{proposition}\label{prop:FI of p form}
If the simple $p-$form $\bar\Omega$ is Frobenius integrable then each form in the sequence
$\ X_1 \righthalfcup \bar\Omega,\ \dots\ ,\ X_{p-1}\righthalfcup\dots\righthalfcup X_1\righthalfcup \bar\Omega$
is Frobenius integrable.
\end{proposition}

 We will now demonstrate a method for locating integrable sub-distributions of non-integrable distribution.

\noindent{\bf{Method for finding the largest integrable sub-distributions of $D$}}\\
Let $D$ be a non-integrable distribution of constant dimension $p$ on $M^n$ and the corresponding co-distribution $D^\perp$ is of dimension $n-p.$ Let $\acute{D}$ be a complementary distribution with $\X(M^n)=D\bigoplus \acute{D}.$ We aim to find the largest integrable sub-distributions of $D$, of dimension $p-\rho$ say, which satisfy the appropriate independence condition.

This is done in stages by using the following steps. We will augment the non-integrable co-distribution $D^\perp$ by adding $\rho$ $1-$forms from $\acute{D}^{\perp}$ and correspondingly reduce the dimension of $D.$ The Frobenius integrability condition for the augmented distribution is divided into two parts, one is  algebraic and the other is differential
 \begin{align}&d\alpha^a\wedge\Omega_\alpha\wedge\Omega_\beta=0, \ \ \ a=1,\ldots, n-p\label{eqn:FI1}\\
&d\beta^b \wedge\Omega_\alpha\wedge\Omega_\beta=0, \ \ \ b =1,\ldots, \rho< p.\label{eqn:FI2}\end{align}
Here $\Omega_{\alpha}:=\alpha^1\wedge\ldots\wedge\alpha^{n-p}$ is a characterizing form for $D^\perp$ and $\beta^b$  are $1-$forms in $\acute{D}^{\perp}$ with $\Omega_\beta:=\beta^1\wedge\dots\wedge\beta^\rho$ and $\bar{\Omega}_V:=\Omega_\alpha\wedge \Omega_\beta\neq 0.$\\

\noindent \textbf{Step 1.}\\
First we add a single $1-$form $\beta\in \acute{D}^{\perp} $ to $D^\perp$ and generate the algebraic conditions \eqref{eqn:FI1} for $\Sp\{\alpha^1,\ldots\alpha^{n-p},\beta\}$ to be Frobenius integrable. If the algebraic conditions can be satisfied then there may be integrable sub-distributions of dimension $p-1$. If the algebraic conditions fail then go to step $2.$\\
\noindent \textbf{Step 2.}\\
 Now add two $1-$forms from $\acute{D}^{\perp}$ to $D^\perp$. The process continues until we have an enlarged co-distribution which may be integrable, if $\rho=p-1$ then the result is trivial. \\
\noindent \textbf{Step 3.}\\
 The next step is to solve the differential conditions \eqref{eqn:FI2}.\\
\noindent\textbf{Step 4.}\label{step3:integration}\\
 Integrate the Frobenius integrable sub-distribution by using  Theorem \ref{sherr result}.\\

It is important to realise that the difficult part in applying the above method is in finding the solution of differential conditions (step $3$). In some cases, simple choice of free parameters (see example \ref{exp:easy phi}) can lead to solutions of PDE, however, in general it is a challenging part of the method. To overcome this difficulty, we  can apply step $3$ and $4$ together, which is called group invariance method to solve the differential conditions on extra $1-$forms. The group invariance method was introduced in \cite[\S 5]{NTGP13} but while working on the current problem we found out that there was an oversight in the original formulation. In section \ref{GIS}, we not only detail the higher dimensional generalisation of group invariance method but also add an extra condition which was missing in \cite[\S 5]{NTGP13}.

\begin{remark}
  The calculations implied by algebraic conditions \eqref{eqn:FI1} and differential conditions \eqref{eqn:FI2} are identical to those implied by Vessiot construction of the involution of maximal order \cite{Vessiot1924} when an independence condition is imposed.
\end{remark}
\begin{remark}
For a given PDE, the dimension of largest integrable, projectible, distributions is the dimension of the space of independent variables $X.$ Then it is straightforward to see how many extra $1-$forms should be added to enlarge the co-distribution which may be integrable.
\end{remark}

\section{Group Invariance Method}\label{GIS}

In this section we discuss the group invariance approach (see \cite{NTGP13} for more details) to solve the differential conditions on extra one forms $\beta^b:$
\begin{equation}\label{eq:diff cond}
  d\beta^b\wedge\Omega_\alpha\wedge\Omega_\beta=0.
\end{equation}

Suppose that we have a general solution of the algebraic conditions \eqref{eqn:FI1} for $\beta^b.$ We now apply step $3$ and $4$ together. Let $\bar{\Omega}_V=\Omega_\alpha\wedge \Omega_\beta$ be a $k-$form, where $k:=n-p+\rho$ and $\Omega_\beta$ comes from the algebraic conditions.  We want $\bar{\Omega}_V$ to be Frobenius integrable. Now suppose there exists a solvable structure of $k$ linearly independent vector fields $X_1,\ldots,X_k$ that we wish to impose on $\bar{\Omega}_V.$ We will impose the condition $$\bar{\Omega}_V(X_1,\ldots,X_k)\neq 0$$ and the conditions that $X_1$ is a symmetry of $\bar{\Omega}_V$, $X_2$ is a symmetry of $X_1\righthalfcup\bar{\Omega}_V$ and so on $X_k$ is a symmetry of $X_{k-1}\righthalfcup\dots\righthalfcup X_1\righthalfcup\bar{\Omega}_V.$ These symmetry conditions are \begin{align}\label{eq:symmcond}
 \mathcal{L}_{X_1}\bar{\Omega}_V &=\ell_1\bar{\Omega}_V,\nonumber\\
 \mathcal{L}_{X_2}(X_1\righthalfcup\bar{\Omega}_V) &=\ell_2(X_1\righthalfcup\bar{\Omega}_V),\\
 \vdots\nonumber\\
 \mathcal{L}_{X_{k}}(X_{k-1}\righthalfcup\dots\righthalfcup X_1\righthalfcup\bar{\Omega}_V) &=\ell_{k}(X_{k-1}\righthalfcup\dots\righthalfcup X_1\righthalfcup\bar{\Omega}_V),\nonumber
 \end{align}
  for some smooth functions $\ell_1,\ldots,\ell_k.$
 
By proposition \ref{prop:FI of p form}, the Frobenius integrability of each form in the sequence $$X_1 \righthalfcup \bar\Omega_V,\ \dots\ ,\ X_{k-1}\righthalfcup \dots\righthalfcup  X_1\righthalfcup  \bar\Omega_V$$is a necessary condition for the Frobenius integrability of $\bar\Omega_V$. And since each form in the sequence has one-form factors in $\Sp\{\alpha^a,\beta^b\}$ their Frobenius integrability  represents a successive simplification of the Frobenius integrability of $\bar\Omega_V$.

We will now further assume that $X_1,\ldots,X_{n-p}$ are linearly independent symmetries of $\Omega_\alpha$ satisfying $X_i\righthalfcup \Omega_\alpha\neq 0,~i=1,\ldots,n-p,$ that is, symmetries of our original PDE. Then we impose the first $n-p$ symmetry conditions \eqref{eq:symmcond} and this puts conditions on $\beta^b$ which simplify the differential conditions \eqref{eq:diff cond}. The remaining $\rho$ symmetries in \eqref{eq:symmcond} of necessity have linearly independent non-zero components in the kernel of $\Omega_\alpha.$  In a similar way, we impose the last $\rho$ symmetry conditions and these conditions will generate more conditions on $\beta^b.$

 If ker $\bar{\Omega}_V\bigoplus \Sp\{X_1,\ldots,X_{k-1}\}$ is Frobenius integrable then the $1-$form
\[\omega:=\frac{X_{k-1}\righthalfcup (X_{k-2}\righthalfcup\ldots\righthalfcup X_{n-p}\righthalfcup\ldots \righthalfcup X_1 \righthalfcup(\Omega_\alpha\wedge\Omega_\beta))}{(\Omega_\alpha\wedge\Omega_\beta)(X_1,\dots,X_k)},\]
is closed by virtue of \eqref{eq:symmcond}. This closure is a necessary condition for the integrability of $\bar\Omega_V.$

If, furthermore, the $\beta^b$ are chosen (see remark \ref{remark:extra 1-forms}) so that $\beta^b(X_i)=0$ for $i=1,\ldots,n-p,~b=1,\ldots,\rho$  then we can simplify this differential condition by writing $\omega$ as
\begin{align*}
 \omega:&=\frac{\Omega_\alpha(X_1,\ldots,X_{n-p})\Omega_\beta(X_{n-p+1},\ldots,X_{k-1})}{\Omega_\alpha(X_1,\ldots,X_{n-p})\Omega_\beta(X_{n-p+1},\dots,X_k)}\\
 &=\frac{\Omega_\beta(X_{n-p+1},\ldots,X_{k-1})}{\Omega_\beta(X_{n-p+1},\dots,X_k)}\in\Sp\{\beta^1,\ldots,\beta^\rho\}.
\end{align*}
This additional requirement on the $\beta^b$ corrects an omission in \cite[\S 5]{NTGP13}.

In forcing the symmetry conditions \eqref{eq:symmcond} and the closure of $\omega$ we have found a closed form in $\Omega_\beta$ and reduced the $\rho$ differential conditions to $\rho-1$ differential conditions. If we have commuting symmetries then we can expect more closed one forms out of this procedure.

In this way we get a solution of the given PDEs which is group invariant. In particular, we obtain $k$ conserved quantities $f^1,\dots,f^k$
\begin{align*}
df^1&=\frac{X_{k-1}\righthalfcup\dots\righthalfcup X_1\righthalfcup \bar{\Omega}_V}{\bar{\Omega}_V(X_1,\dots,X_k)}=\frac{\Omega_\beta(X_{n-p+1},\ldots,X_{k-1})}{\Omega_\beta(X_{n-p+1},\dots,X_k)},\\
df^2&\equiv\frac{X_{k}\righthalfcup\dots\righthalfcup X_1\righthalfcup \bar{\Omega}_V}{\bar{\Omega}_V(X_1,\dots,X_k)}\quad \text{mod}\ df^1,\\
&\vdots\qquad\qquad\vdots\\
df^k&\equiv\frac{X_{k}\righthalfcup\dots\righthalfcup X_2\righthalfcup \bar{\Omega}_V}{\bar{\Omega}_V(X_1,\dots,X_k)} \quad \text{mod}\ df^1,\dots,df^{k-1}.
\end{align*}
\begin{remark}\label{remark:extra 1-forms}
  We can always choose such one-forms $\beta^b$ which satisfy the condition $\beta^b(X_i)=0$ for $i=1,\ldots,n-p,~b=1,\ldots,\rho$ by adding appropriate multiple of $\alpha^a.$

  For example, set $\bar{\beta}^b=\beta^b+r^b_a\alpha^a.$ The condition $\bar{\beta}^b(X_i)=0,$ implies
  \begin{align*}
    0=\bar{\beta}^b(X_i)=\beta^b(X_i)+r^b_a\alpha^a(X_i).
  \end{align*}

  This system always has a solution for $r^b_a$ because the $X_i$ are linearly independent and $X_i\righthalfcup\Omega_\alpha\neq 0.$
\end{remark}
\section{Multi-dimensional second order PDEs}
\noindent In this section, we examine second order PDEs in higher dimensions. Firstly, we apply the above method to a class of second order PDEs of one dependent variable and three independent variables. Secondly, we apply the same technique to a second order PDEs of one dependent variable and four independent variables. Finally, some examples are given for the system of second order PDEs.
\subsection{${\bf{3-}}$Dimensional second order PDEs}
Consider a partial differential equation of the form
\begin{equation}\label{3GPDE}
  G(x,y,z,u,u_x,u_y,u_z,u_{xx},u_{xy},u_{xz},u_{yy},u_{yz},u_{zz})=0,
\end{equation}
in three independent variables $x,~y,~z$ and one dependent variable $u.$
The embedded submanifold
\[S:=\{(x,y,\ldots,u_{zz})\in J^2(\R^3,\R)~ |~ u_{zz}-F(x,y,z,u,u_x,u_y,u_z,u_{xx},u_{xy},u_{xz},u_{yy},u_{yz})=0\}\]
is a subset of $J^2(\R^3,\R).$ A local solution of the PDE is a $12-$dimensional locus of $J^2(\R^3,\R)$ described by the map $i:S~\hookrightarrow~J^2(\R^3,\R),$ i.e. \[i:(x,y,z,u,u_x,u_y,u_z,u_{xx},u_{xy},u_{xz},u_{yy},u_{yz})\hookrightarrow (x,y,z,u,u_x,u_y,u_z,u_{xx},u_{xy},u_{xz},u_{yy},u_{yz},F).\]
We can solve \eqref{3GPDE} for any of the second order derivatives. We consider the equation of the following form
\begin{equation}\label{3PDE}
  u_{zz}=F(x,y,z,u,u_x,u_y,u_z,u_{xx},u_{xy},u_{xz},u_{yy},u_{yz}).
\end{equation}
We study the solutions of the given PDE by studying the integral submanifolds $N$ of the pulled-back contact system $D^{\perp}_V$
\begin{align}
\theta^1&:=du-u_xdx-u_ydy-u_zdz,\nonumber\\
 \theta^2&:=du_x-u_{xx}dx-u_{xy}dy-u_{xz}dz,\label{3contact dist} \\
 \theta^3&:=du_y-u_{xy}d x-u_{yy}dy-u_{yz}dz \nonumber\\
 \theta^4&:=du_z-u_{xz}d x-u_{yz}dy-Fdz \nonumber,
 \end{align}
which project down to $X\subset \R^3.$ If $N$ satisfies the independence condition $dx\wedge dy\wedge dz| _N \neq 0$ and has a tangent space that annihilates the distribution $D^{\perp}_V$, then $i(N) \subset  J^2(\R^3,\R)$ is the $2-$graph of a solution of \eqref{3PDE}. The corresponding dual distribution $D_V$  is generated by
\begin{align}
V_1&:=\vp{x}+u_{x}\vp{u}+u_{{x}{x}}\vp{u_{x}}+u_{xy}\vp{u_{y}}+u_{xz}\vp{u_{z}},\nonumber\\
V_2&:=\vp{y}+u_{y}\vp{u}+u_{xy}\vp{u_{x}}+u_{yy}\vp{u_{y}}+u_{yz}\vp{u_{z}},\label{3Vessiot dist}\\
V_3&:=\vp{z}+u_{z}\vp{u}+u_{xz}\vp{u_{x}}+u_{yz}\vp{u_{y}}+F\vp{u_{z}},\nonumber\\
V_4&:=\vp{u_{{x}{x}}},~V_5:=\vp{u_{xy}},~V_6:=\vp{u_{yy}},\nonumber\\
V_7&:=\vp{u_{xz}},~V_8:=\vp{u_{yz}}.\nonumber
\end{align}
The Vessiot distribution is not Frobenius integrable. To find the largest integrable sub-distributions we apply the method of section $3$.

For this class of PDEs we need to add five $1-$forms $\phi^1,\ldots,\phi^5\in \acute{D}^{\bot}$ in $D^{\perp}_V$ because the reduced Vessiot distribution $D_{red}$ must have dimension $3$ in order to project to the tangent distribution of a $2-$graph. Without loss of generality, we assume that
\begin{align}\label{eqn:phi1,2}
\phi^1&:=du_{xx}-a_{1} dx-a_{2}dy-a_{3}dz,\nonumber\\
\phi^2&:=du_{xy}-a_{4}dx-a_{5}dy-a_{6}dz,\nonumber\\
\phi^3&:=du_{yy}-a_{7}dx-a_{8}dy-a_{9}dz,\\
\phi^4&:=du_{xz}-a_{10}dx-a_{11}dy-a_{12}dz,\nonumber\\
\phi^5&:=du_{yz}-a_{13}dx-a_{14}dy-a_{15}dz.\nonumber
\end{align}
These particular forms of the $\phi$'s guarantee the independence condition.  Now solve the algebraic conditions
\[d\theta^a\wedge\Omega_\theta\wedge\Omega_\phi=0,\quad a=1,\ldots,4.\]
Specifically,
\begin{align}
  a_2&=a_4,~~a_3=a_{10},~~a_6=a_{11}=a_{13},~~a_5=a_7,~~a_9=a_{14},\nonumber\\
  a_{12}&=u_xF_u+u_{xx}F_{u_x}+a_1F_{u_{xx}}+a_4F_{u_{xy}}+a_{10}F_{u_{xz}}+u_{xy}F_{u_y}+a_{7}F_{u_{yy}}\nonumber\\&+a_{13}F_{u_{yz}}+u_{xz}F_{u_{z}}+F_x,\label{3PDE:alg sol}\\
 a_{15}&=u_yF_u+u_{xy}F_{u_x}+a_2F_{u_{xx}}+a_5F_{u_{xy}}+a_{11}F_{u_{xz}}+u_{yy}F_{u_y}+a_{8}F_{u_{yy}}\nonumber\\&+a_{14}F_{u_{yz}}+u_{yz}F_{u_{z}}+F_y.\nonumber
\end{align}
At this stage we need an expression for $F.$ For this purpose we consider the Boyer Finley equation. This equation has been studied by many authors within different frameworks \cite{Dunajski2001,Sergey,boyereq,Martina2001,ward1990}.
\begin{example}
Consider the Boyer Finley equation
$u_{zz}=e^{u_{xx}}-u_{yy}.$
From \eqref{3PDE:alg sol}, we have
\begin{align*}
 a_2&=a_4,~~a_3=a_{10},~~a_6=a_{11}=a_{13},~~a_5=a_7,\\
  a_9&=a_{14},~~a_{12}=a_1e^{u_{xx}}-a_7,~ a_{15}=a_2e^{u_{xx}}-a_8.
\end{align*}

After solving the algebraic conditions, extra $1-$forms are
\begin{align*}
  &\phi^1=du_{xx}-a_1dx-a_2dy-a_3dz,\quad \phi^2=du_{xy}-a_2dx-a_5dy-a_6dz,\\
  &\phi^3=du_{yy}-a_5dx-a_8dy-a_9dz\quad \phi^4:=du_{xz}-a_3dx-a_6dy-(a_1e^{u_{xx}}-a_5)dz,\\
  & \phi^5:=du_{yz}-a_6dx-a_9dy-(a_2e^{u_{xx}}-a_8)dz.
\end{align*}
At this stage we have seven unknown functions $a_1,~a_2,~a_3,~a_5,~a_6,~a_8$ and $a_9$ of $12$ variables. We apply the group invariance method to solve the differential conditions.

Here, $\bar{\Omega}_V$ is a $9-$form so we need nine linearly independent symmetries to impose the symmetry conditions \eqref{eq:symmcond}. The symmetries of the original PDE are
\begin{align*}
X_1&:=\vp{u},~~ X_2:=y\vp{u}+\vp{u_y},~~ X_3:=-z\vp{u}-\vp{u_z},~~ X_4:=x\vp{u}+\vp{u_x},~~ X_5:=\vp{x},\\
X_6&:=x\vp{u_y}+\vp{u_{xy}}+y\vp{u_x}+xy\vp{u},~~ X_7:=-(x\vp{u_z}+\vp{u_{xz}}+z\vp{u_x}+xz\vp{u}),\\
X_8&:=-2z\vp{u_z}+2\vp{u_{yy}}+2y\vp{u_{y}}+(y^2-z^2)\vp{u},~~X_9:=-(y\vp{u_z}+\vp{u_{yz}}+z\vp{u_y}+yz\vp{u}).
\end{align*}
The linearly independence condition on symmetries implies \[\bar{\Omega}_V(X_1,\ldots,X_{9})=-2a_1\neq 0.\] Also, these symmetries satisfy the condition $X_i\righthalfcup \Omega_\theta\neq 0$ and $\phi^b(X_i)= 0,~i=1,\ldots,4,~b=1,\ldots,5.$
The first symmetry condition $\mathcal{L}_{X_1}\bar{\Omega}_V =\ell_1\bar{\Omega}_V$ implies
\begin{align*}
&\pd{a_1}{u}=0,~\pd{a_2}{u}=0,~\pd{a_3}{u}=0,~\pd{a_5}{u}=0,~\pd{a_6}{u}=0,~\pd{a_8}{u}=0,~\pd{a_9}{u}=0.
\end{align*}
In a similar fashion, we apply the next eight symmetry conditions and these conditions imply that $a_1,~a_2,~a_3,~a_5,~a_6,~a_8$ and $a_9$ are function of three variables $y,~z,u_{xx}.$ The last one form should be closed by virtue of symmetry conditions, i.e.
\begin{align*}
  \omega^{9}:=\frac{\Omega_\phi(X_5,X_6,X_7,X_8)}{\Omega_\phi(X_5,X_6,X_7,X_8,X_9)} .
\end{align*}
The closure condition $d\omega^{9}=0$  and Frobenius integrability of other one-forms will generate further conditions on unknown functions. One possible solution of this over-determined system is
\[a_1=b_1e^{-u_{xx}},~a_2=0,a_3=0,a_5=(y+z)b_2+b_3,~a_6=(y-z)b_2+b_4~a_8=\frac{b_2}{b_1}e^{u_{xx}},~a_9=\frac{b_2}{b_1}e^{u_{xx}}+b_5,\]
where $b_1$ is non zero constant and $b_2,\ldots,b_5$ are arbitrary constants. For the sake of convenience, we choose these constants to be equal to one.

Note that, $\phi^1,\ldots,\phi^5$ satisfy the differential conditions and the transverse condition. The corresponding $D_{red}$ is generated by
\begin{align*}
\bar{V}_1&:=V_1+(y+z+1)V_5+e^{u_{xx}}V_6+(y-z+1)V_7+(e^{u_{xx}}+1)V_8,\\
\bar{V}_2&:=V_2+(y+z+1)V_5+e^{u_{xx}}V_6+(y-z+1)V_7+(e^{u_{xx}}+1)V_8,\\
\bar{V}_3&:=V_3+(y-z+1)V_5+(e^{u_{xx}}+1)V_6-(y+z)V_7-e^{u_{xx}}V_8.
\end{align*}
After integration we obtain the following invariant functions:
\begin{align*}
f^1&:=u-xu_x+xyu_{xy}+xzu_{xz}-yu_y+yzu_{yz}-zu_{z}-xy^2z+xyz^2-xyz+xu_{xx}e^{u_{xx}}\\
&+\frac{1}{12}[2z^3-6y^2z+4xz^3-6xy^2-4xy^3-9x^2-6z^2u_{yy}+6y^2u_{yy}+e^{u_{xx}}(6z^2+2z^3+6yz^2\\
&-6y^2z-2y^3+6x)-6u_{xx}e^{2u_{xx}}],\\
f^2&:=xu_{xy}-u_{y}-yu_{yy}+zu_{yz}-xy-xyz-xz-yz-yze^{u_{xx}}+\frac{1}{2}((z^2-y^2)e^{u_{xx}}-xy^2+xz^2),\\
f^3&:=xu_{xz}-u_{z}-zu_{yy}+yu_{yz}-xy+xyz+(yz+z)e^{u_{xx}}+\frac{1}{2}((z^2-y^2)e^{u_{xx}}-xy^2+xz^2-y^2-z^2),\\
f^4&:=yz^2-yz-y^2z+yu_{xy}-u_x+zu_{xz}+(u_{xx}-1)e^{u_{xx}}+\frac{1}{6}(2z^3-3y^2-2y^3),\\
f^5&:=x-e^{u_{xx}},~~f^6:=z+y+yz-u_{xy}+\frac{1}{2}(y^2-z^2),~~f^7:=y-yz-u_{xz}+\frac{1}{2}(y^2-z^2),\\
f^8&:=\frac{1}{2}(z-u_{yy}+(y+z+1)e^{u_{xx}})~~f^9:=y-u_{yz}+(y-z+1)e^{u_{xx}}.
\end{align*}
\noindent The lifted solution on $S$ is a common level set of these functions:
$$\left\{ p\in S:\ f^\alpha(p)=c_\alpha\right\}.$$
This projects to
 \begin{align*}
  u&=c_1-c_2y-zc_3-xc_4-c_5yz-c_6xy-c_7xz-c_8(y^2-z^2)-c_9yz+xyz+\frac{1}{12}[2xy^3+6xy^2z\\&+6xy^2-6xyz^2-2xz^3+6y^2z-2z^3-9x^2+c_5(2z^3+6yz^2-6y^2-6y^2z-2y^3+18x)\\
  &+(6c_5^2-12xc_5+6x^2)\log(x-c_5)].
\end{align*}
This is a group invariant solution and the reader can check the invariance.
\end{example}

\subsection{${\bf{4-}}$Dimensional second order PDEs}
Consider a partial differential equation of the form
\begin{equation}
  G(t,x,y,z,u,u_t,u_x,u_y,u_z,u_{tx},u_{ty},u_{tz},u_{tt},u_{xx},u_{xy},u_{xz},u_{yy},u_{yz},u_{zz})=0,
\end{equation}
in four independent variables $t,~x,~y,~z$ and one dependent variable $u.$
The embedded submanifold
\[S:=\{(t,x,\ldots,u_{yz},u_{zz})\in J^2(\R^4,\R) | u_{zz}-F(t,x,y,z,u,u_t,\ldots,u_{yz})=0\}\]
is a subset of $J^2(\R^4,\R).$ A local solution of the PDE is a $18-$dimensional locus of $J^2(\R^4,\R)$ described by the map $i:S\hookrightarrow J^2(\R^4,\R),$ i.e. \begin{align*}
  i&:(t,x,y,z,u,u_t,u_x,u_y,u_z,u_{tx},u_{ty},u_{tz},u_{tt},u_{xx},u_{xy},u_{xz},u_{yy},u_{yz})\\ \hookrightarrow & (t,x,y,z,u,u_t,u_x,u_y,u_z,u_{tx},u_{ty},u_{tz},u_{tt},u_{xx},u_{xy},u_{xz},u_{yy},u_{yz},F).
\end{align*}
The restriction (or pullback by inclusion) of the contact system $D^{\perp}_V$ on $J^2(\R^3,\R)$ is given by
\begin{align}
\theta^1&:=du-u_tdt-u_xdx-u_ydy-u_zdz, \quad\quad \theta^2:=du_t-u_{tt}dt-u_{tx}dx-u_{ty}dy-u_{tz}dz,\nonumber\\
 \theta^3&:=du_x-u_{tx}dt-u_{xx}dx-u_{xy}dy-u_{xz}dz,\quad \theta^4:=du_y-u_{ty}dt-u_{xy}d x-u_{yy}dy-u_{yz}dz\label{4contact dist} \\
 \theta^5&:=du_z-u_{tz}dt-u_{xz}d x-u_{yz}dy-Fdz. \nonumber
 \end{align}
If the integral submanifolds $N$ of the restricted contact system satisfies the independence condition $dt\wedge dx\wedge dy\wedge dz| _N \neq 0$ and has a tangent space that annihilates the distribution $D^{\perp}_V$, then $i(N) \subset  J^2(\R^4,\R)$ is the $2-$graph on a solution of the given PDE. The corresponding dual distribution $D_V$ is generated by
\begin{align}
V_1&:=\vp{t}+u_{t}\vp{u}+u_{tt}\vp{u_{t}}+u_{{t}{x}}\vp{u_{x}}+u_{ty}\vp{u_{y}}+u_{tz}\vp{u_{z}},\nonumber\\
V_2&:=\vp{x}+u_{x}\vp{u}+u_{tx}\vp{u_{t}}+u_{{x}{x}}\vp{u_{x}}+u_{xy}\vp{u_{y}}+u_{xz}\vp{u_{z}},\nonumber\\
V_3&:=\vp{y}+u_{y}\vp{u}+u_{ty}\vp{u_{t}}+u_{xy}\vp{u_{x}}+u_{yy}\vp{u_{y}}+u_{yz}\vp{u_{z}},\label{4Vessiot dist}\\
V_4&:=\vp{z}+u_{z}\vp{u}+u_{tz}\vp{u_{t}}+u_{xz}\vp{u_{x}}+u_{yz}\vp{u_{y}}+F\vp{u_{z}},\nonumber\\
V_5&:=\vp{u_{xx}},~V_6:=\vp{u_{xy}},~V_7:=\vp{u_{yy}},~V_8:=\vp{u_{xz}},V_9:=\vp{u_{yz}},\nonumber\\
V_{10}&:=\vp{u_{tx}},~V_{11}:=\vp{u_{ty}},~V_{12}:=\vp{u_{tz}},~V_{13}:=\vp{u_{tt}}.\nonumber
\end{align}
This distribution is not Frobenius integrable, so we now apply the method (section $3$) of finding integrable sub-distributions.
For this class of PDEs we need to add nine $1-$forms $\phi^1,\ldots,\phi^9\in \acute{D}^{\bot}$ in $D^{\perp}_V$ because the reduced Vessiot distribution must have dimension $4$ in order to project to the tangent distribution of a $2-$graph. Once again we demonstrate with an example:
\begin{example}\label{exp:easy phi}

Consider the so-called second heavenly equation \cite{Sergey}
$u_{tx}+u_{yz}+u_{xx}u_{zz}-u_{xz}^2=0.$

\noindent The Vessiot distribution and co-distribution is given by equations \eqref{4contact dist} and \eqref{4Vessiot dist}  with $F=\frac{1}{u_{xx}}(u_{xz}^2-u_{yz}-u_{tx}).$

For instance, if we choose $\phi^1:=du_{xx},~ \phi^2:=du_{xy},~ \phi^3:=du_{yy},~ \phi^4:=du_{xz},~ \phi^5:=du_{yz},~ \phi^6:=du_{tx},~ \phi^7:=du_{ty},~ \phi^8:=du_{tz},~ \phi^9:=du_{tt}$, we obtain a Frobenius integrable sub-distribution.

The symmetries of reduced Vessiot distribution $D_\text{red}$ are
\begin{align*}
X_1&:=\vp{u},\quad X_2:=\vp{z}, \quad X_3:=-y\vp{u}-\vp{u_y},\quad X_4:=-x\vp{u}-\vp{u_x},\\
X_5&:=-t\vp{u}-\vp{u_t},\quad X_6:=-(x\vp{u_y}+\vp{u_{xy}}+y\vp{u_x}+xy\vp{u}),\\
X_7&:=y\vp{u_{z}}+\vp{u_{yz}}+z\vp{u_{y}}-t\vp{u_x}-\vp{u_{tx}}-x\vp{u_t}+(yz-xt)\vp{u},\\
X_8&:=2\vp{u_{yy}}+2y\vp{u_y}-y^2\vp{u},\quad X_9:=-(t\vp{u_y}+\vp{u_{ty}}+y\vp{u_t}+ty\vp{u}),
\end{align*}
\begin{align*}
X_{10}&:=2\vp{u_{tt}}+2t\vp{u_t}-t^2\vp{u},\\ X_{11}&:=-u_{z}\vp{u_{z}}+u_{yy}\vp{u_{yy}}+u_{xy}\vp{u_{xy}}+u_{xx}\vp{u_{xx}}-u_{tz}\vp{u_{tz}}-u_{tt}\vp{u_{tt}}
-u_t\vp{u_t}\\&-u\vp{u}-y\vp{y}-x\vp{x},\\
X_{12}&:=x\vp{u_{z}}+2u_{xz}\vp{u_{xz}}+4u_{xy}\vp{u_{yy}}+2u_x\vp{u_{y}}+\vp{u_{xz}}+2u_{xx}\vp{u_{xy}}+z\vp{u_{x}}\\
&+2u_{tx}\vp{u_{ty}}-xz\vp{u}-2y\vp{x},\\
X_{13}&:=\frac{1}{2u_{xx}}[-u_{xx}u_{xz}\vp{u_{yz}}-2u_{xx}u_{xy}\vp{u_{yy}}-u_xu_{xx}\vp{u_{y}}-u_{xx}^2\vp{u_{xy}}+(u_{xz}^2-u_{yz}-u_{tx})\vp{u_{tz}}\\
&+(u_{xx}u_{yz}-u_{xx}u_{tx})\vp{u_{ty}}+u_{xx}u_{xz}\vp{u_{tx}}+2u_{xx}u_{tz}\vp{u_{tt}}+u_zu_{xx}\vp{u_t}-tu_{xx}\vp{z}+yu_{xx}\vp{x}],\\
X_{14}&:=u_{z}\vp{u_{z}}-2u_{yz}\vp{u_{yz}}-4u_{yy}\vp{u_{yy}}-u_y\vp{u_{y}}-u_{xz}\vp{u_{xz}}-3u_{xy}\vp{u_{xy}}-2u_{xx}\vp{u_{xx}}\\
&-u_{tz}\vp{u_{tz}}-3u_{ty}\vp{u_{ty}}-2u_{tx}\vp{u_{tx}}-2u_{tt}\vp{u_{tt}}+2u\vp{u}+2x\vp{x}+2t\vp{t}+3y\vp{y}+z\vp{z}.
\end{align*}
and after integration we obtain the following invariant functions:
\begin{align*}
f^1&:=tu_t+xu_x+yu_y+zu_z-u-xyu_{xy}-xzu_{xz}-yzu_{yz}-tzu_{tz}-tyu_{ty}-txu_{tx}\\
&-\frac{1}{2}(t^2u_{tt}+y^2u_{yy})+\frac{1}{2u_{xx}}(z^2u_{yz}^2-z^2u_{tx}-z^2u_{xz}^2-x^2u_{xx}^2),\\
f^2&:=z+\frac{1}{u_{tz}}(tu_{tt}+xu_{tx}+yu_{ty}-u_t),\quad f^3:=u_y-tu_{ty}-xu_{xy}-zu_{yz}-yu_{yy},\\
f^4&:=zu_{xz}+tu_{tx}-u_x+yu_{xy}+xu_{xx},\\
f^5&:=u_t-tu_{tt}-xu_{tx}-yu_{ty}+\frac{u_{xx}u_{tz}}{(u_{xz}^2-u_{tx}-u_{yz})}(xu_{xz}-u_z+yu_{yz}+tu_{tz}),\\
f^6&:=-u_{xy},\quad f^7:=u_{tx},\quad f^8:=\frac{1}{2}u_{yy},\quad f^9:=u_{ty},\quad f^{10}:=\frac{1}{2}u_{tt},\\
f^{11}&:=\log (u_{tx}+u_{yz})-\log {u_{xx}},~f^{12}:=\frac{1}{2}\log(2u_{xz}+1),~f^{13}:=\frac{2u_{xx}u_{tz}}{u_{yz}-u_{xz}^2+u_{tx}},\\
f^{14}&:=-\frac{1}{2}\log (u_{tx}+u_{yz}).
\end{align*}
The common level sets project to
\begin{align*}
  u&=c_3y-c_1-c_4x-c_6xy+c_7(tx-yz)+c_{10}t^2+c_8y^2+c_9 t-2c_5c_{13}^{-1}z+e^{-2c_{14}}yz\\
  &+\frac{1}{2}[e^{c_{11}}(c_{13}tz-c_2c_{13}t+2c_2z-z^2)+e^{2c_{12}}xz+e^{-c_{11}-2c_{14}}x^2-xz]\\
  &+\frac{1}{8}[e^{c_{11}+2c_{14}}(z^2-2c_2z-c_{13}tz+c_{2}c_{13}t)+e^{c_{11}+2c_{12}+2c_{14}}(4c_2z-2z^2+2c_{13}tz-2c_{2}c_{13}t)
  \\&+e^{c_{11}+4c_{12}+2c_{14}}(z^2-2c_2z-c_{13}tz+c_{2}c_{13}t)].
  \end{align*}
\end{example}
\newpage
\subsection{Systems of PDEs}
We apply the group invariance method to a system of nonlinear second order PDEs of two dependent variables and two independent variables.
\begin{example}
Consider a system of two nonlinear second order PDEs \cite{VAM2009}
\begin{align}
 3u_{xx}u_{yy}^3+1=0,\quad  v_{xx}-\frac{v_{yy}}{u_{yy}^4}=0.
 \end{align}
 The pulled-back  contact system $D^{\perp}_V$ on $J^2(\R^2,\R^2)$ with $F_1=-\frac{1}{3u_{yy}^3}$ and $F_2=\frac{v_{yy}}{u_{yy}^4}, ~u_{yy}\neq 0$ is generated by
\begin{align}
\theta^1&:=du-u_xdx-u_ydy,\quad\theta^2:=du_v-v_{x}dx-v_{y}dy,\nonumber\\
 \theta^3&:=du_x-F_1dx-u_{xy}dy,\quad \theta^4:=du_y-u_{xy}d x-u_{yy}dy\label{syscontact dist} \\
 \theta^5&:=dv_x-F_2dx-v_{xy}dy,\quad\theta^6:=dv_y-v_{xy}d x-v_{yy}dy.\nonumber
 \end{align}
It is easy to see that we need to add extra four $1-$forms in $D^{\perp}_V$ and the solution of algebraic conditions \eqref{eqn:FI1} implies
\[a_{3}  =a_{2},\quad a_7=a_6,\quad a_4=a_1u_{yy}^4,\quad a_8= \frac{1}{u_{yy}}(4a_4v_{yy}+a_5u_{yy}^5).\]
After solving the algebraic conditions, extra $1-$forms are
\begin{align*}
  &\phi^1=du_{xy}-a_1dx-a_2dy,\quad \phi^2=du_{yy}-a_2dx-a_1u_{yy}^4dy,\\
  &\phi^3=dv_{xy}-a_5dx-a_6dy, \quad \phi^4:=dv_{yy}-a_6dx-(4a_1v_{yy}u_{yy}^3-a_5u_{yy}^4)dy.
\end{align*}
At this stage we have four unknown functions $a_1,~a_2,~a_5$ and $a_6$ of $10$ variables. We impose the symmetry conditions \eqref{eq:symmcond} on $\bar{\Omega}_V.$ The symmetries of the original PDE are
\begin{align*}
X_1&:=\vp{u},\quad X_2:=\vp{v},\quad X_3:=-(\vp{u_x}+x\vp{u}),\quad X_4:=-(\vp{u_y}+y\vp{u}),\\
X_5&:=-(\vp{v_x}+y\vp{v}),\quad X_6:=-(\vp{v_y}+y\vp{v}),~X_7:=-(x\vp{u_y}+\vp{u_{xy}}+y\vp{u_x}+xy\vp{u}),\\
X_8&:=-(x\vp{v_y}+\vp{v_{xy}}+y\vp{v_x}+xy\vp{v}),\quad  X_9:=-(v_{yy}\vp{v_{yy}}+v_{xy}\vp{v_{xy}}+v_y\vp{v_y}+v_x\vp{v_x}+v\vp{v}),\\
X_{10}&:=-4v_{yy}\vp{v_{yy}}-2v_{xy}\vp{v_{xy}}-u_{yy}\vp{u_{yy}}+u_{xy}\vp{u_{xy}}-2v_{y}\vp{v_{y}}+3u_x\vp{u_{x}}+u_y\vp{u_y}+3u\vp{u}+2y\vp{y}.
\end{align*}
The linearly independence condition on symmetries implies \[\bar{\Omega}_V(X_1,\ldots,X_{10})=u_{yy}v_{yy}(2a_1yu_{yy}^3-1)\neq 0.\]
The first eight symmetry conditions \eqref{eq:symmcond} imply that the four unknown functions depend only on three variables $y,u_{yy}$ and $v_{yy}.$ The last two symmetry conditions implies
\begin{align*}
&\pd{a_2}{v_{yy}}=0,~\pd{a_1}{v_{yy}}=0,~v_{yy}\pd{a_6}{v_{yy}}+a_6=0,~v_{yy}\pd{a_5}{v_{yy}}+a_5=0,\\
&u_{yy}\pd{a_2}{u_{yy}}-2y\pd{a_2}{y}-a_2=0,~u_{yy}\pd{a_1}{u_{yy}}-2y\pd{a_1}{y}-a_1=0.
\end{align*}
 Now impose the closure of the last one form, which generates further conditions on unknown functions. The last one form is
\begin{align*}
  \omega^{10}:=\frac{\Omega_\phi(X_7,X_8,X_9)}{\Omega_\phi(X_7,X_8,X_9,X_{10})}.
\end{align*}
Note that, $X_9$ and $X_{10}$ are commuting symmetries and we can expect to get two closed one forms. By imposing the closure conditions and remaining two Frobenius integrability conditions and solving this over-determined system, we have
\begin{align*}
  &a_1=0,~a_2=0,~a_5=0,~a_6=0,\quad \text{or}\\
  &a_1=\frac{b^1}{yu_{yy}^3},~a_2=\frac{u_{yy}}{x+b^2},~a_5=0,~a_6=\frac{v_{yy}}{x+b^3},
\end{align*}
where $b^1=0$ or $b^1=\frac{-1}{3}$ and $b^2,~b^3$ are arbitrary constants.

Choosing $a_1=0,~a_2=0,~a_5=0,~a_6=0,$ after integration, we obtain the following conserved quantities:
\begin{align*}
f^1&:=xu_x+yu_y-xyu_{xy}-u+\frac{1}{6u_{yy}^3}(x^2-3y^2u_{yy}^4),\\
f^2&:=v-xv_x-yv_y+xyv_{xy}+\frac{v_{yy}}{2u_{yy}^4}(x^2+y^2u_{yy}^4),\\
f^3&:=u_x-yu_{xy}+\frac{x}{3u_{yy}^3},\quad  f^4:=xu_{xy}-u_y+yu_{yy},\\
 f^5&:=v_{x}-yv_{xy}-\frac{xv_{yy}}{u_{yy}^4},\quad f^6:=xv_{xy}-v_y+yv_{yy},\\
 f^7&:=u_{xy},~ f^8:=-v_{xy},~f^9:=\log v_{yy}-4\log u_{yy},~f^{10}:=-\log u_{yy}.
\end{align*}
The common level sets project to
\begin{align*}
  u&=c_3x-c_4y-c_1+c_7xy+\frac{1}{6}(3y^2e^{-c_{10}}-x^2e^{3c_{10}}),\\
  v&=c_2-c_6y+c_5x-c_8xy+\frac{1}{2}(x^2e^{c_{9}}+y^2e^{c_{9}-4c_{10}}).
\end{align*}

\end{example}
Now we present a group invariant solution of a system of nonlinear second order PDEs of one dependent variable and three independent variables.
\begin{example}

Consider a system of three nonlinear second order PDEs
\begin{align}
 u_{xx}-u_{zz}(u_{yy}u_{zz}+u_{yz}^2)=0,\quad  u_{xy}+u_{yy}u_{zz}+\frac{1}{2}u_{yz}^2=0,\quad  u_{xz}+u_{yz}u_{zz}=0.
\end{align}
The pulled-back contact system $D^{\perp}_V$ on $J^2(\R^3,\R)$ is generated by
\begin{align}
\theta^1&:=du-u_xdx-u_ydy-u_zdz,\quad \theta^2:=du_x-F_1dx-F_2dy-F_3dz,\nonumber\\
 \theta^3&:=du_y-F_2d x-u_{yy}dy-u_{yz}dz,\quad \theta^4:=du_z-F_3d x-u_{yz}dy-u_{zz}dz,\label{syscontact dist}
 \end{align}
where $F_1=u_{zz}(u_{yy}u_{zz}+u_{yz}^2),~F_2=-u_{yy}u_{zz}-\frac{1}{2}u_{yz}^2$ and $F_3=-u_{yz}u_{zz}.$  By adding three extra $1-$forms  in $D^{\perp}_V$ and solving algebraic conditions \eqref{eqn:FI1}, we have
\begin{align*}
&a_{6}  =a_{8},\quad a_5=a_3,\quad a_9=0,\quad a_7=-a_6u_{zz},\\
&a_4=-a_5u_{zz}-a_8u_{yz},\quad a_1= -a_2u_{zz}-a_5u_{yz}-a_8u_{yy}.
\end{align*}
After solving the algebraic conditions, extra $1-$forms are
\begin{align*}
  &\phi^1=du_{yy}+(a_2u_{zz}+a_3u_{yz}+a_6u_{yy})dx-a_2dy-a_3dz,\\
  &\phi^2=du_{yz}+(a_3u_{zz}+a_6u_{yz)}dx-a_3dy-a_6dz,~~\phi^3=du_{zz}+a_6u_{zz}dx-a_6dy.
\end{align*}
At this stage we have three unknown functions $a_2,~a_3$ and $a_6$ of $10$ variables. We impose the symmetry conditions \eqref{eq:symmcond} on $\bar{\Omega}_V.$ The symmetries of the original PDE are
\begin{align*}
X_1&:=\vp{u},\quad X_2:=-\vp{u_y}-y\vp{u},\quad X_3:=-\vp{u_x}-x\vp{u},\\
X_4&:=-\vp{u_z}-z\vp{u},\quad X_5:=\vp{z},\quad X_6:=\vp{x},\quad X_7:=\vp{y}.
\end{align*}
The linearly independence condition on symmetries implies $\bar{\Omega}_V(X_1,\ldots,X_7)=-a_6^3u_{yy}\neq 0$ and the symmetry conditions imply that $a_2,~a_3$ and $a_6$ are functions of three variables $u_{yy},~u_{yz}$ and $u_{zz}.$

Similarly, we can impose the closure of the last one-form which generates more conditions on unknown functions.
By solving these conditions with Frobenius integrability conditions, we have
\[a_2=0,~a_3=0,~a_6=b_1u_{yy},\] where $b_1$ is a non zero constant.

After integration, we obtain the following conserved quantities:
\begin{align*}
f^1&:=u-xu_x-yu_y-zu_z+\frac{1}{2(xb_1u_{yy}-1)}(2xyu_{yy}u_{zz}-x^2u_{zz}^2u_{yy}-y^2u_{yy}\\
&+yz^2b_1u_{yy}-x^2u_{zz}u_{yz}^2+xyu_{yz}^2+2xzu_{yz}u_{zz}-2yzu_{yz}-z^2u_{zz}),\\
f^2&:=u_{y}-\frac{1}{b_1}u_{zz}-\frac{u_{yz}^2}{2b_1u_{yy}},\quad f^3:=-u_x-\frac{1}{2b_1}u_{zz}^2-\frac{u_{zz}u_{yz}^2}{2b_1u_{yy}},~f^4:=u_z-\frac{u_{yz}u_{zz}}{b_1u_{yy}},\\
f^5&:=z-\frac{u_{yz}}{b_1u_{yy}},\quad f^6:=-x+\frac{1}{ b_1u_{yy}},\quad f^7:=y-\frac{u_{zz}}{b_1u_{yy}}.
\end{align*}
The common level sets project to
 \begin{align*}
 u&=c_1+c_2y-c_3x+c_4z+ \frac{1}{2c_6b_1(c_6 + x)}(c_5^2 c_6b_1y + c_5^2 c_7b_1x + 2c_5c_6c_7b_1z - 2c_5c_6b_1yz - 2c_6c_7y\\
 &- c_6c_7b_1z^2+ c_6y^2  + c_6b_1yz^2 - c_7^2 x).
\end{align*}
\end{example}
This is a group invariant solution.
\begin{remark}
We remark that, in the context of hyperbolic second order PDEs in the plane, finding the extra (closed) 1-forms in our method is equivalent
to finding the second order Darboux invariants. In other words, the calculations one goes through to construct the Darboux invariants agrees with our method and we refer readers to \cite{Goursat} where all this and more can be seen in great detail.
\end{remark}

\section{Discussion and future work}
We have demonstrated in this paper a method for the solution of second order partial differential equations of two and more than two independent variables. Vessiot theory is used to generate the integrable sub-distributions which satisfy the appropriate independence condition. By this method, we can find a finite parameter subset of an infinite-parameter family of the integrable sub-distributions. In this process, we subdivide the Frobenius condition into two parts, one being algebraic and one differential. Solving the former is straight forward, the most challenging task is to satisfy the differential conditions (step $3$). After solving the algebraic part, various parameters remains, which must satisfy the differential conditions. For some PDEs, simple choices of the free parameters can lead to solutions. Then a solvable structure of each integrable sub-distribution is used to integrate the reduced Vessiot distribution, giving a local solution of the original system as the parameterised integral submanifold of the sub-distribution. But, in general, the differential conditions cannot be solved explicitly. Vessiot has given an existence theory \cite{Vessiot1924} of such Frobenius integrable distributions in the analytic category. We have used solvable symmetry structure of the original PDEs to explicitly obtain the integrable sub-distributions of Vessiot distributions, giving the group invariant solution.  Then a natural question arises ``what is the relationship between the construction of group invariant solutions using this method and classical group invariant solutions?" The answer to this question is given in  \cite[\S VI]{sherr}.

It is worth pointing out again the utility of our geometric approach in the investigation of some natural questions concerning conserved quantities and symmetries. Having solved the problem of finding the maximal integrable sub-distributions of the Vessiot distribution, we can now ask, for example, does my PDE admit constant curvature solutions?\\
The answer is: at least locally if and only if $dK\wedge df^1\wedge\ldots\wedge df^p\equiv0,$  where $K$ is a Gaussian curvature of the graph of a solution of the given PDE.
But more usefully it does not admit a constant curvature solution if and only if
$dK\wedge \Omega_{\theta}\wedge \Omega_{\phi}\neq 0,$ where $\Omega_\phi$ only has to satisfy the algebraic conditions.
Moreover, if $K$ is not constant on the graph then the one dimensional integral manifolds of $dK\wedge \Omega_{\theta}\wedge \Omega_{\phi}$ are the curves on the graph along which $K$ is constant and their local existence is guaranteed.

 The future tasks involve the generalisation to higher order systems. Also, we will work on the construction of symmetries leaving boundary and initial value problems invariant and the construction of the corresponding invariant solutions of PDEs. We believe that introducing initial and boundary conditions in this framework will introduce some conditions on free parameters after solving the algebraic conditions. We also believe that if we impose a solvable symmetry structure which respects initial/boundary conditions, then the resultant group invariant solutions will also respect the conditions. Moreover, it would be interesting to classify (in invariant terms) those PDEs for which a solvable structure can be imposed in order to determine Frobenius integrable sub-distributions.
\section*{Acknowledgements}
The author gratefully acknowledges the support of a La Trobe University  postgraduate research award and the kind hospitality of Department of Mathematics, Aarhus University, Denmark. She would also like to thank Geoff Prince and the anonymous referee for  many valuable suggestions which have improved the presentation of this paper significantly.


\providecommand{\bysame}{\leavevmode\hbox to3em{\hrulefill}\thinspace}
\providecommand{\MR}{\relax\ifhmode\unskip\space\fi MR }
\providecommand{\MRhref}[2]{%
  \href{http://www.ams.org/mathscinet-getitem?mr=#1}{#2}
}
\providecommand{\href}[2]{#2}

\end{document}